\documentclass[a4paper,12pt]{amsart}

\usepackage{amsmath,amsthm,amssymb}

\font\Bbb=msbm10 at 10pt

\newcommand{\eps}{{\varepsilon}}

\newcommand{\ZZ}{{\hbox{\Bbb{Z}}}}
\newcommand{\hh}{{\rm h}}

\newcommand{\beq}[1]{\begin{equation}\label{#1}}
\newcommand{\eeq}{\end{equation}}
\newcommand{\vv}{{\rm v}\kern .3pt}

\newtheorem{theorem}{Theorem}[section]
\newtheorem{lemma}[theorem]{Lemma}
\newtheorem{proposition}[theorem]{Proposition}
\newtheorem{cor}[theorem]{Corollary}
\newtheorem{conj}[theorem]{Conjecture}
\newcommand{\bep}{\begin{proof}}
\newcommand{\eep}{\end{proof}}
\newcommand{\bec}[1]{\begin{cor}\label{#1}}
\newcommand{\eec}{\end{cor}}
\newcommand{\bepr}[1]{\begin{proposition}\label{#1}}
\newcommand{\eepr}{\end{proposition}}
\newcommand{\bel}[1]{\begin{lemma}\label{#1}}
\newcommand{\eel}{\end{lemma}}
\newcommand{\bet}[1]{\begin{theorem}\label{#1}}
\newcommand{\eet}{\end{theorem}}
\newcommand{\beco}[1]{\begin{conj}\label{#1}}
\newcommand{\eeco}{\end{conj}}

\def\card{{\rm Card}\kern .3pt}
\def\hh{{\rm h}\kern .3pt}

\def\bQ{{\bf Q}}
\def\bR{{\bf R}}

\title{On a  conjecture on exponential Diophantine equations}

\author{Mihai Cipu \and Maurice Mignotte}\thanks{
The first author has been partially supported by the CEEX Program 
of the Romanian Ministry of Education, Research and Youth, Grant 
2-CEx06-11-20/2006.}

\date{}

\address{Institute of Mathematics\\
Romanian Academy  \\
   P.O Box 1--764\\
RO-014700 Bucharest\\
Romania and 
Universit\'{e} Louis Pasteur\\
          U. F. R. de Math\'{e}matiques\\
          7, rue Ren\'{e} Descartes\\
          67084 Strasbourg Cedex\\
          France
}

\email{mihai.cipu@imar.ro}

\address{Universit\'{e} Louis Pasteur\\
          U. F. R. de Math\'{e}matiques\\
          7, rue Ren\'{e} Descartes\\
          67084 Strasbourg Cedex\\
          France}
\email{mignotte@math.u-strasbg.fr}

\begin{document}

\begin{abstract}
We deal with a conjecture of Terai (1994)  asserting that if
$a$, $b$, $c$ are fixed coprime integers with $\min (a,b,c)>1$ such that
$a^2+b^2=c^r$ for a certain odd integer $r>1$, then the equation
$a^x+b^y=c^z$ has only one solution in positive integers with $\min
(x,y,z)>1$.
   Co-operation man-machine is needed for the proof.
\end{abstract}

\maketitle

\section{The problem}\label{sint}

Let $a$, $b$, $c$ be fixed coprime integers with $\min (a,b,c)>1$.
In 1933, Mahler~\cite{ma} developed a $p$-adic equivalent of the
Thue-Siegel method to prove that the equation
\beq{ec1}
a^x+b^y=c^z
\eeq
has finitely many solutions $(x,y,z)$ in positive integers. His
method is ineffective in the sense that it gives no indication on
the number of possible solutions for a fixed triple $(a,b,c)$.
Such an information has been obtained only in particular instances.
Thus,  Sierpi\'nski~\cite{si} showed that  $(2,2,2)$ is the
unique solution in positive integers to the equation $3^x+4^y=5^z$.
In the same journal, Je\'smanowicz~\cite{je} conjectured the
unicity of the solution to Eq.~(\ref{ec1}) in case $(a, b, c)$
is a Pythagorean triple. This conjecture is still open, despite
the efforts of many authors.

In analogy to Je\'smanowicz's conjecture, Terai~\cite{tet} stated
that Eq.~(\ref{ec1}) always has at most one solution in positive 
integers. Simple examples disproving this statement have been
found by Cao~\cite{cao}, who attempted to remedy the situation by 
adding the hypothesis $\max (a,b,c) >7$. It turns out that this 
condition is not sufficient to entail the thought-for unicity. A 
family of counterexamples have been pointed out by Le~\cite{le},
who also stated the following variant of Terai's conjecture.

\begin{conj}\label{conj1}
For given coprime integers $a$, $b$, $c>1$, the Diophantine
equations~\emph{(\ref{ec1})} has at most one  solution in
integers $x$, $y$, $z>1$.
\end{conj}

Much work has been devoted to the case when Eq.~(\ref{ec1}) has a 
solution of the form $(2,2,r)$, with $r$ greater than 1 and odd. 
This implies in particular that $c$ is odd and exactly one of $a$, 
$b$ is even. For the sake of definiteness, suppose that $a$ is 
even and therefore $b$ is odd. Most of the recent results  concerns 
the case $a\equiv 2 \pmod 4$,  $b\equiv 3 \pmod 4$. The conjecture 
is established in this case under one of the
following additional hypotheses:
\begin{enumerate}
\item[$(\alpha)$] (Terai~\cite{ter2}) $b\equiv 3 \pmod 8$, $b\ge 30a$
and $\left(\frac{a}{d}\right)=-1$, where  $d>1$ is a divisor of $b$ 
and $\left(\frac{a}{d}\right)$ denotes the Jacobi symbol,

\item[$(\beta)$] (Cao~\cite{cao}) $c$ is a prime power,

\item[$(\gamma)$] (Cao-Dong~\cite{dcao}) $b\ge 25.1 a$,

\item[$(\delta)$] (Le~\cite{le}) $c>3\cdot 10^{27}$ and $r>7200$.
\end{enumerate}

\noindent Further partial confirmations of the conjecture are 
referred to in the papers cited above.

Contrary to what it is claimed in~\cite{le}, the last result 
quoted above does not imply that the conjecture holds with the 
exception of finitely many pairs $(c,r)$.
One of the aims of this paper is to prove that indeed there are 
at most finitely many values for which the conjecture can be refuted.
On the way we shall prove other results for the positive solutions
to the Diophantine simultaneous equations
\beq{ec2}
a^2+b^2=c^r, \qquad a^2+b^y=c^z,
\eeq
where
\beq{ec2bis}
r  , z  \ {}>1\ \mathrm{are~ odd}, \ a \equiv 2
\pmod 4,  \   b \equiv 3\pmod 4, \ \mathrm{and} \ \gcd (a,b)=1.
\eeq

As a consequence of our deliberations, improvements on the results
$(\alpha )$--$(\delta )$ are obtained. Our proofs approach these 
cases from a different perspective and are much shorter than the 
published ones, although they involve a harder computational 
component. We give here a rough description of our procedure. In 
the hypotheses of our work, $c$ is a sum of two coprime squares. 
We generate all such decompositions for $c$ up to $4\cdot10^{10}$ 
with the help of Cornacchia's algorithm (see~\cite{ni} and~\cite{ba} 
for very simple proofs of its correctness). We notice that, when 
compared to the obvious method (for $c$ fixed and $1\le u<\sqrt c$ 
test whether $c-u^2$ is a square), for our range of values 
Cornacchia's algorithm is more than ten times faster.

This description is vague; details are given in the third
section, after we recall classical facts, some of them going
back at least  to Lagrange. Additional information on the
putative solutions of the Diophantine system~(\ref{ec2}) are
given  in Section~\ref{sdoi}. Section~\ref{strei} contains
the proofs of our main results, among which are the following.

\bet{tea}
If the Diophantine equation $X^x+Y^y=Z^z$  has a solution with
$X=a\equiv 2\pmod 4$,  $Y= b \equiv 3\pmod 4$, $Z=c$, $x=2$, $y=2$
and $z=r$ odd, where $\gcd (a,b)=1$, then this is the only 
solution in positive integers, with the possible exception of 
finitely many values $(c,r)$.
\eet

\bet{tec}
If  $a$ or $b$  is a prime power then the
system~\emph{(\ref{ec2})}  has   no solutions subject
to restrictions  from~\emph{(\ref{ec2bis})}.
\eet
The last part of the paper is devoted to improvements of
bounds on the parameters associated to a putative solution
to system~(\ref{ec2}). They are meant to  shrink the search
domains for the components of a solution to a manageable
size according to the present-day technology.

Although Terai's conjecture remains open, we have pushed the
analysis further than ever before; and there is significant
hope that our results can be improved by either  complementing
them with brand new ideas or dedicated computations.

\section{Arithmetic restrictions}\label{sunu}

We use a result of Lagrange (1741), \emph{Le\c cons sur le calcul 
des fonctions}, which makes recurrent appearance in the study of 
Diophantine equations, as well as in the theory of finite fields, 
Chebyshev polynomials and many other areas of mathematics. For the 
sake of completeness, we sketch its proof.

\bel{le3}
Let $X$ and $Y$ two commuting indeterminates and let $n\ge 1$ be a
positive integer. Then
\[
X^n+Y^n = \sum_{j=0}^{\lfloor n/2\rfloor} c_{n,j}\, (-XY)^j 
(X+Y)^{n-2j},
\]
where the $c_{n,j}$ are nonnegative integers which are defined
recursively  by
\[
c_{n,j}=0 \quad \mathit{if}\quad j<0 \quad \mathit{or}\quad
j>\lfloor n/2\rfloor,
\]
\[
c_{1,0}=1, \ c_{2,0}=1, \ c_{2,1}=2, \quad \mathit{and}\quad
c_{n+1,j}=c_{n,j}+c_{n-1,j-1}\  \mathit{for} \ n\ge 2.
\]
More precisely,
\[
c_{n,j}=\frac{(n-j-1)! \, n}{(n-2j)!\, j!}.
\]
\eel \bep
The result is obvious for $n\le2$, including the initial values
$c_{1,0}=1$, $c_{2,0}=1$, $c_{2,1}=2$.
The general case can be obtained by induction from the formula
\[
X^{n+1}+Y^{n+1}=(X^n+Y^n)(X+Y)-XY(X^{n-1}+Y^{n-1}),
\]
which implies the recursive relation
\[
c_{n+1,j}=c_{n,j}+c_{n-1,j-1}\  {\rm for}\ n\ge 2.
\]
By completely working out the details, one can get the closed form 
for the coefficients $c_{n,j}$.
\eep

In the present situation, using Corollary~\ref{cormord} below, we 
get   expressions for $a$, $b$ and~$b ^{y/2}$.

\bec{cor1}
   The values of $a$ and $b$ satisfy
\[
a = \pm u \sum_{j=0}^{(r-1)/2} c_{r,j}\, (-c)^j (4v^2)^{(r-1)/2-j},
\quad
b = \pm v \sum_{j=0}^{(r-1)/2} c_{r,j}\, c^j (-4v^2)^{(r-1)/2-j}
\]
and
\[
a = \pm u_1 \sum_{j=0}^{(z-1)/2} c_{z,j}\, (-c)^j (4v_1^2)^{(r-1)/2-j},
\quad
b ^{y/2}= \pm v_1 \sum_{j=0}^{(z-1)/2}\, c_{z,j} c^j
(-4v_1^2)^{(z-1)/2-j}.
\]
\eec

  From the last formula it follows that
\[
b ^{y/2}\equiv \pm v_1 z c^{(z-1)/2} \pmod{v_1^3},
\]
in particular
\[
p \mid v_1 \ \Longrightarrow  \  \vv_p(v_1)\ge y/2\ge 3
\quad {\rm if}\ \gcd(p,z)=1.
\]

\section{Bounds for $a$, $b$, $c$}\label{sdoi}

From our standard hypotheses on  $a$, $b$, $c$ stated in 
Eq.~(\ref{ec2}) and~(\ref{ec2bis})  it follows that 
$c\equiv 5\pmod8$. As it is well-known (see, e.g., 
Lemma~\ref{lemord} below), the first equation from~(\ref{ec2}) 
implies that there exist positive integers $u$ and $v$ such that
\beq{ec3bis}
c=u^2+v^2.
\eeq
By~\cite{cao},  we also may suppose that $c$ has at least two 
prime divisors. Then  it is easily seen  that one has $c\ge 85$.

Other useful facts are given by the next result, proved in 
several places, for instance in~\cite{le}.

\bel{le0}
With the above notation and hypotheses, let $(x,y,z)$ be a 
solution to~\emph{(\ref{ec1})} with $(x,y,z)\ne (2,2,r)$. Then 
$x=2$, $y\equiv 2\pmod 4$, $y\ge 6$ and $z$ is odd.
\eel

We shall repeatedly use the well-known structure of integers 
satisfying the first  equation from~(\ref{ec2}).

\bel{lemord}
If $X$, $Y$ and $Z$ are coprime positive integers such that
\[
X^2+Y^2=Z^n,
\]
where $n$ is an odd integer and $X$ is even, then there exist 
coprime positive integers $u$ and $v$, with $u$ even and $v$ odd, 
and $\lambda_1$, $\lambda_2\in \{-1,1\}$ such that
\[
X+Y\sqrt{-1} = \lambda_1 \eps^n, \quad \eps=u+v\lambda_2\sqrt{-1},
\quad Z=u^2+v^2.
\]
Moreover, if $\varepsilon=|\varepsilon|\,e^{\theta\sqrt{-1}/2}$ 
then
\[
X=Z^{n/2}|\cos(n\theta/2)|, \quad
Y=Z^{n/2}|\sin(n\theta/2)|.
\]
\eel

The former part is proved as in Mordell's 
book~\cite[pp.122--123]{mor}; the later assertion is an  obvious 
consequence of the preceding formulas.

\medskip

In the present case, using the relations $a^2+b^2=c^r$ and 
$a^2+b^y=c^z$,  we get exponential expressions for $a$, $b$ 
and~$b ^{y/2}$.

\bec{cormord}
There are positive integers $u$, $v$,  $u_1$, $v_1$, with $u$, 
$u_1$  even and  $v$, $v_1$ odd,  such that $c=u^2+v^2=u_1^2+v_1^2$  
and
\[
a=\frac{1}{ 2} |\eps ^r +\bar \eps^r| = c^{r/2}|\cos(r \xi)|,
\quad
b= \frac{1}{ 2} |\eps ^r -\bar \eps^r|= c^{r/2}|\sin(r \xi)|
\]
and
\[
a= \frac{1}{ 2}
  |\eps_1 ^z +\bar \eps_1^z| = c^{z/2}|\cos(z \xi_1)|,
\quad
b^{y/2} =  \frac{1}{ 2} |\eps_1 ^z -\bar \eps_1^z| = 
c^{z/2}|\sin(z \xi_1)|,
\]
where $ \eps=u+v\sqrt{-1}$, $ \eps_1=u_1+v_1\sqrt{-1}$,
$\tan \xi =v/u$, and $\tan \xi_1 =v_1/u_1$.
\eec

\medskip

Let us come back to the notation of the above lemma and put  
$\eps=e^{i\theta/2}$ and
\[
\alpha:=\varepsilon/\bar \varepsilon=\frac{ u^2-v^2+2 
uv\,  \sqrt{-1}}{ u^2+v^2} \, e^{\theta\, \sqrt{-1}}.
\]
Since $\alpha$ is a root of the irreducible integer polynomial
\[
(u^2+v^2) T^2-2(u^2-v^2)T+u^2+v^2,
\]
whose Mahler's measure is equal to $u^2+v^2$,  the absolute
logarithmic height of $\alpha$ is
\[
\hh(\alpha)= \frac{1}{ 2}\log (u^2+v^2)=  \frac{1}{ 2}
\log Z.
\]
We  also have
\[
\min\{X,Y\}\ge \frac{Z^{n/2}}{ \pi}\, \min_{k'\in \ZZ}
|n\theta - k'\pi|.
\]
Let $k$ be an integer such that $\,\min_{k'\in \ZZ}|
n\theta - k'\pi| =  |n\theta - k\pi|$ and put
\[
\Lambda = n\log \alpha - k \log (-1).
\]
Then
\[
\min\{X,Y\}\ge \frac{Z^{n/2}}{\pi}\,|\Lambda|
\quad \hbox{and }\quad
\min\{X,Y\}\ge 0{.}99\,{Z^{n/2}} \,\min\{|\Lambda|,0{.}001\},
\]
where $\Lambda$ is a linear form in two logarithms of algebraic 
numbers.

\subsection{A first application of linear forms}\label{ss1}

In a number field $K$ embedded in the complex field, containing 
a root   of unity $\zeta=e^{i\pi/m}$, where $m$ is maximal, 
and a  number $\alpha$ of  modulus one
which is not a root of unity, a linear form
\[
\Lambda = n\log \alpha - i k \pi
\]
as above can be written as
\[
\Lambda = n\log \alpha -  m k \log \zeta.
\]
We remark that changing $\alpha$ into a suitable 
$\alpha \zeta^\ell$ if necessary we can assume that 
$|\log\alpha|\le \pi/(2m)$. We may work under this hypothesis 
without changing the notation because
$\hh(\zeta ^\ell\alpha)=\hh(\alpha)$.

On using the main
result of Laurent-Mignotte-Nesterenko~\cite{lmn},
it is possible to prove that
\beq{ec3}
z>55000 \quad{\rm implies}\quad a>c^{z/(2\sqrt 3)}.
\eeq

On using  relations $a^2+b^y=c^z$ and $c=u_1^2+v_1^2$, by a 
computation
of a suitable continued fraction we verify that
\beq{ec4}
{\rm for} \ 85\le c <4\cdot 10^{10},\ z>10  \quad{\rm implies}
\quad a> c^{z/(2\sqrt 3)}.
\eeq
Similarly, from $a^2+b^2=c^r$ and $c=u^2+v^2$ for some integers
$u$, $v$ which may be different from $u_1$, $v_1$, we obtain
\beq{ec4bis}
{\rm for} \ 85\le c <4\cdot 10^{10},\ r>10  \quad{\rm implies}
\quad b> c^{r/(2\sqrt 3)}.
\eeq

This information is exploited in conjunction with the following
remarks.

\bel{le0bis}
Assume both conditions~\emph{(\ref{ec2})} and~\emph{(\ref{ec2bis})}
are fulfilled. Then:

a) If for some $\mu >0$ one has  $a\ge c^{z/\mu}$ then $2z< \mu r$.

b) If for some $\lambda >0$ one has  $b\ge c^{r/\lambda}$ then 
$yr<  \lambda z$.

c) If $\mu _1>0$, $\mu _2>0$ are such that $\mu _1\mu _2 \le 2y$
then  $a\ge c^{z/\mu _1}$ and $b\ge c^{r/\mu _2}$ cannot
simultaneously hold. In particular,
\[
a<c^{z/(2\sqrt 3)}  \quad or \quad b< c^{r/(2\sqrt 3)}.
\]
\eel
  \bep
$a)$ From $a\ge c^{z/\mu}$ and $a^2+b^2=c^r$ it readily follows 
that $c^{2z/\mu}<c^{r}$.

$b)$ If  $b\ge c^{r/\lambda}$ then $c^{yr/\lambda}\le b^y<c^z$.

$c)$ The first assertion follows directly from $a)$ and $b)$.
The last part follows from this because $y\ge 6$ by
Lemma~\ref{le0}.
\eep

Using this lemma, we could rule out the small values of $r$
and $z$ (precisely, those with $2< r< z <10$)
and prove that $c$ cannot be comparatively small. After around 
two weeks of computation we could verify the following result.

\bel{le0bibi}
Assume the system of equations~\emph{(\ref{ec2})} has solutions satisfying~\emph{(\ref{ec2bis})}. Then $ c>4\cdot 10^{10}$.
\eel

\subsection{A second application of linear forms}\label{ss2}

From now on we consider  $c>4\cdot 10^{10}$ without further 
explicitly mentioning it. In order to improve the bounds on $r$ 
and $z$ obtained in the previous subsection, we apply a very 
recent result of Laurent~\cite{lau}.

\bel{laur}
Consider a nonzero linear form
\[
\Lambda = b_1\log \alpha_1 - b_2\log \alpha_2,
\]
where $\alpha_1$ and $\alpha_2$ are nonzero algebraic numbers,
both different from $1$, and $b_1$ and $b_2$ are positive
integers. Put
\[
D =[\bQ(\alpha_1,\alpha_2):\bQ]/[\bR(\alpha_1,\alpha_2):\bR].
\]
Let $K$ be an integer $\ge 3$, $L$ an integer $\ge 2$,
$R_1$, $R_2$, $S_1$, $S_2$ positive integers.
Let $\rho $ and $\mu$ be  real
numbers with $\rho >1$ and $1/3\le \mu\le1$. Put
$\,R=R_1+R_2-1$, $S=S_1+S_2-1$, $N=KL$,
\begin{align}
g &= \frac{1}{4}-\frac{N}{ 12RS}, \qquad
\sigma=\frac{1+2\mu-\mu^2}{ 2}, \notag
\\
b &= \frac{\bigl((R-1)b_2+(S-1)b_1\bigr)}{ 2}
\left(\prod_{k=1}^{K-1}k!\right)^{-2/(K^2-K)}. \notag
\end{align}
Let $a_1$, $a_2$ be positive real numbers such that
\[
a_i \ge \rho \,|\log \alpha _i| - \log |\alpha _i|+2D\,
{\rm h} \kern  .3pt(\alpha_i) ,
\]
for $i=1$, $2$. Suppose that:
\[
\card \bigl\{ \alpha _1^r\,\alpha _2^s \,;\,
0\le r<R_1, \,0\le s<S_1 \bigr\} \ge L,
\leqno (I)
\]
\[
\card \bigl\{ rb_2+sb_1 \,;\, 0\le r<R_2,\,0\le s<S_2\bigr\}
> (K-1)L
\leqno (II)
\]
and also that
\[
K(\sigma L-1)\log\rho-
(D+1)\log N -D(K-1)\log b -gL\,(R a_1+S
a_2)\,>\,c(N),  \leqno (III)
\]
where
\[
c(N)=\frac{2}{ N}\,\log\left(N!\, N^{-N+1} \bigl(e^N +(e-1)^N
\bigr)\right).
\]
Then
\[
\vert\Lambda'\vert\ge\rho^{-\mu KL },
\]
where
\[
\Lambda^\prime=\Lambda \cdot\max
\left\{\frac{LSe^{LS\vert\Lambda\vert/(2b_2)}}{2b_2}, \,
\frac{LRe^{LR\vert\Lambda\vert/(2b_1)}}{2b_1}\right\}.
\]
\eel

In our case $\alpha_1=\alpha$ (up to a power of $\sqrt{-1}\,$),  
$\alpha_2=\sqrt{-1}$, $b_1=r$ or $z$, and $b_2=k$. 
(To work with the  linear form  associated to the  relation
$a^2+b^y=c^z$ we only need to take above $b_1=z$ instead of
$b_1=r$.)
For $c=4\cdot 10^{10}+5$, we choose the parameters as follows: $L=8$,
$\rho=7{.}7$, $\mu=0{.}56$, $K = \lceil m L a_1a_2 \rceil$,
$ R_1=4$,  $S_1=2$,  $R_2= \lceil \sqrt m L a_2\rceil $,
and  $S_2=\lceil  (1+(K-1)L)/R_2\rceil $,
where $m=0{.}1166$, and we get
\[
|\Lambda |> c^{-0{.}2113 r }\quad {\rm for}\ r\ge 771,
\]
which implies
\[
  a > c^{  z/( 2\sqrt 3)}\quad {\rm and}\
   b > c^{  r/( 2\sqrt 3)}\quad {\rm for}\ r\ge 771.
\]
Taking into account Lemma \ref{le0bis}$b)$, one concludes that
$r\le 769$.

Now, combining Lemma \ref{le0bis}$a)$ and Lemma \ref{le0bis}$c)$, 
we   see that if the system has a  solution then $r\le 769$ and 
$z\le 983$. The detailed   argument is the following: we apply 
Laurent's result twice, a first computation for $z\ge 985$ gives 
an upper bound for $\mu_2$ which combined with   part $a)$ implies 
$r\ge 641$; then a second computation for  $r\ge 641$ gives  an
 upper bound for $\mu_1$ with $\mu_1\mu_2<12$, and  
part $c)$ leads to a contradiction. Thus $z\le 983$. 
Moreover, it is easy to check that the greater $c$,  
the better our estimates, so that
the conclusion holds for all $c>4\cdot 10^{10}$.

Arguing
in the same way, we can establish tighter bounds for $r$ and
$z$, provided a higher lower bound on $y$ is available.

\bel{le0tri}
If the Diophantine system~\emph{(\ref{ec2})} has solutions
satisfying~\emph{(\ref{ec2bis})}  then in all cases
\[
r\le 769\quad  and  \quad z\le 983.
\]
Moreover
\[
y \ge 10 \ \Longrightarrow  \ r\le 539\quad and \quad  z\le 759,
\]
\[
y \ge 14 \ \Longrightarrow  \ r\le 461 \quad and \quad  z\le 681,
\]
\[
y \ge 18 \ \Longrightarrow  \ r\le 419 \quad and \quad  z\le 647,
\]
\[
y \ge 22 \ \Longrightarrow  \ r\le 395 \quad and \quad  z\le 627,
\]
and
\[
y \ge 602 \ \Longrightarrow  \ r\le 263  \quad and \quad z\le 539.
\]
\eel

\subsection{Elementary lower bounds on $b$}\label{ss4}

Let $\eps=u+iv=|\eps|e^{i\xi}$, where $c=u^2+v^2$,
with $u$ even,  and
$|\eps|=\sqrt c$. Then $\tan  \xi =v/u$ and
\[
b=\frac{1}{2}| \eps^r-\bar\eps^r|=c^{r/2}|\sin (r\xi)|,
\]
with $r\ge 3$.
In this subsection we derive lower bounds on $b$ from
lower bounds on $v$.

\bel{le1}
With the above notation, one has
\[
r\le \pi/\xi -1 \ \Longrightarrow \ b\ge v\, c^{(r-1)/2}\ge v\,c.
\]
In particular, one gets $b\ge v\, c^{(r-1)/2}$ whenever
$r\le u\pi/v -1$.
\eel \bep
The hypothesis $3\le r\le \pi/\xi -1$ implies that
$\xi\le \pi /4$ and
$3\xi\le r \theta\le\pi -\xi$,
and therefore $\sin (r\xi) \ge \sin \xi =v/\sqrt c$.

For the last part, note that the hypothesis $r\le \pi/\xi -1$
holds if  $r\le u\pi/v -1$ because $0<\xi <\tan \xi =v/u$.
\eep

Despite its innocuous appearance, lemma just proved plays an
important role in subsequent reasonings. Thus, $v\le 925$
implies $u/v> 216$ (recall our standing hypothesis 
$c>4\cdot 10^{10}$) and then the previous lemma gives 
$b\ge v \, c^{(r-1)/2} \ge c^{r/3}$ (since  $3\le r$).
Having in view  Lemma~\ref{le3}, it follows that one always 
has $b\ge v$. Therefore,  $b\ge 925$.

More importantly, with the help of Lemma~\ref{le1}  we shall 
derive a strikingly sharp bound for the quotient $y/z$.

\bel{le1bis}
We always have
\[
b \ge \frac{\pi}{r+1} \left(1+\frac{\pi^2}{(r+1)^2}\right)^{-1/2} 
\!\sqrt c
\]
and
\[
y < z \left(2+\frac{9{.}982}{ \log b} \right).
\]
Moreover, if $y>600$ then
\[
y < z \left(2+\frac{8{.}863}{ \log b} \right).
\]
In particular, it always holds
\[
y< 1778.
\]
\eel

\bep
 From our previous study we know that
\beq{ecmnb}
b \ge \begin{cases}
c^{(r-1)/2}\ge c &  \hbox{\rm if \ $(r+1)v<\pi u$},
\\
v &  {\rm otherwise}.
\end{cases}
\eeq
Notice that $(r+1)v\ge \pi u$ implies
\[
c \le \left(1+\frac{(r+1)^2}{\pi^2}\right)v^2,
\]
so that in all cases $b$ satisfies
\[
b \ge \frac{\pi}{r+1} \left(1+\frac{\pi^2}{(r+1)^2}\right)^{-1/2}
\! \sqrt c.
\]
Now we consider the upper bounds for $y$. From \ref{ecmnb} we get
\[
c \ \le \begin{cases}
b &  \hbox{\rm if \ $(r+1)v<\pi u$},
\\
v^2 + \frac{(r+1)^2}{\pi^2} v^2 &  {\rm otherwise}.
\end{cases}
\]
Hence,
\[ c\le \left(1+ \frac{(r+1)^2}{\pi^2}\right)b^2.
\]
Using the inequality $b^y<c^z$ one gets
\beq{ec9}
y < z \left(2+\frac{\log  \left(1+(r+1)^2/\pi^2 \right)}{ \log b} 
\right).
\eeq
If $y\le 10$ the second estimate of the lemma  is trivial,  
hence we suppose $y\ge 14$.
Then $r<462$, and after a simple computation we get the 
stated inequality.

When $y$ is greater than $600$ we know from  Lemma~\ref{le0tri}
that $r\le 263$ and the third estimate follows. The last one
is deduced by using the fact that $b$ is at least $925$ and $z$
is less than $540$ whenever $y$ is at least $600$.
\eep

\subsection{Estimates on $a$}\label{ss7}
Our next goal is to obtain some estimates on $a$. Put $b=c^\lambda$.
The information we have up to know allows us to conclude that
$1/2-(\log 1800)/\log c<\lambda < r/2$. We use this knowledge to 
prove the following.

\bel{le2}
Put $a=b^{(y-\lambda')/2}$. Then $\lambda'$ is positive and 
satisfies
\[
   \lambda' > \frac{\log c}{ \log b}\left(z-r-10^{-22}  \right) 
>  \frac{2}{ r}\left(2-10^{-22}  \right).
\]
\eel \bep
 From the second equation in~(\ref{ec2}) we get
\[
c^z(1-c^{-z+r})<b^y<c^z
\]
and since $z\ge r+2$ this implies
\[
z \log c + \log(1-c^{-2})<y\log b<z\log c,
\]
while the first equation in~(\ref{ec2}) and the definition of 
$\lambda'$ imply
\[
(y-\lambda')\log b < r \log c.
\]
Hence
\[
z \log c -10^{-21} <  r \log c + \lambda' \log b ,
\]
and therefore
\[
2 \le z-r < 10^{-22} + \lambda' \,\frac{\log b}{ \log c}10^{-22}  
+ \lambda'  \lambda,
\]
by the definition of $\lambda$. In other words
\[
   \lambda' > \frac{1}{ \lambda}\left(z-r-10^{-22}  \right),
\]
and in particular
\[
   \lambda' > \frac{2}{ r}\left(2-10^{-22}  \right)>0.
\]
\eep

\section{Main results}\label{strei}

Recall the result of Corollary~\ref{lemord}: we have seen that  
$c=u^2+v^2=u_1^2+v_1^2$ for some positive
integers, with $u$, $u_1$ even and  $v$, $v_1$ odd, and that
\[
a=\frac{1}{ 2} |\eps ^r +\bar \eps^r| =  \frac{1}{ 2}
  |\eps_1 ^z +\bar \eps_1^z|,
\quad
b= \frac{1}{ 2} |\eps ^r -\bar \eps^r|,
\quad
b^{y/2} =  \frac{1}{ 2} |\eps_1 ^z -\bar \eps_1^z|,
\]
where $ \eps=u+v\sqrt{-1}$ and $\eps_1=u_1+v_1\sqrt{-1}$. 
It follows that, up to a sign,
$a$, $b$ and $b^{y/2}$ are values of binary linear recursive
sequences. If $(u,v)=(u_1,v_1)$ then the term $b^{y/2}$ has no
primitive divisors, so that on checking tables of binary Lucas
sequences having terms without primitive divisors given
in~\cite{bhv} and~\cite{abid}  we recover Cao's result~\cite{cao}
mentioned in Introduction.

\bet{leBHV}
If  $c$  is a prime power then the system~\emph{(\ref{ec2})}
has   no solutions subject to restrictions
from~\emph{(\ref{ec2bis})}.
\eet

\medskip

Now we are in a position to prove that the conjecture holds
perhaps with the exception of finitely many pairs $(c,r)$.

Subtracting the two equations from~(\ref{ec2}) results in
the Diophantine equation
\beq{ec12}
b^y- b^2 =c^z -c^r.
\eeq
Since $6\le y$ and $5\le z$, for fixed exponents $(y,r,z)$
one  gets an algebraic curve of positive genus.
The absolute irreducibility and the genus of the curve
defined by Eq.~(\ref{ec12}) are given by a theorem of
Davenport,  Lewis and  Schinzel~\cite{dls}.

\bel{leo2}
Let $f(X)$ and $g(Y)$ be polynomials with integral
coefficients of degree $n>1$ and respectively $m>1$.
Let $D(\lambda )= \mathrm{disc}(f(X)+\lambda )$ and
$E(\lambda )=\mathrm{disc}(g(Y)+\lambda )$. Suppose
there are at least $n/2$ distinct roots of $D(\lambda )=0$
for which $E(\lambda )\ne 0$. Then $f(X)-g(Y)$ is
irreducible over the complex field. Further, the genus of
the curve $f(x)-g(y)=0$ is positive except possibly  when
$m=2$ or $m=n=3$. Apart from these possible exceptions,
the  equation $f(x)-g(y)=0$ has at most finitely many
integral solutions.
\eel

Stickelberger's formula~\cite{sti} (cf.~\cite{sw}) for 
the discriminant of a trinomial gives
\[\mathrm{disc}(b^y- b^2+\lambda ) = -\lambda \left(
y^{y/2}\lambda^{y/2-1} -2(y-2)^{y/2-1}\right)^2,
\]
\[\mathrm{disc}(c^z -c^r+\lambda )= (-1)^{z(z-1)/2}
\lambda^{r-1}\left( z^z\lambda^{z-r}-(z-r)^{z-r}r^r \right),
\]
so that the last quoted result applies.

Combining these classical facts with some of our results
in the previous sections, we obtain the main result of
the paper.
\bet{teo1}
If the Diophantine equation $X^x+Y^y=Z^z$  has a solution with
$X=a\equiv 2\pmod 4$,  $Y= b \equiv 3\pmod 4$, $Z=c$, $x=2$, $y=2$
and $z=r$ odd, where $\gcd (a,b)=1$, then this is the only solution
in positive integers, with the possible exception of finitely many
values $(c,r)$.
\eet \bep
For each fixed pair of odd numbers $(r,z)$, $1<r<z$, any
solution to the system~(\ref{ec2}) subject to~(\ref{ec2bis})
corresponds to an integer point on a curve~(\ref{ec12})
of positive genus. By Siegel's seminal paper~\cite{sie},
such an equation has only finitely many integral solutions.
According to Lemma~\ref{le0tri}, in any compatible
system~(\ref{ec2}) one has $r<770$ and $z\le 983$.
Moreover, $y$ is bounded from above by $1800$ (see 
Lemma~\ref{le1bis}).
Therefore, a compatible system~(\ref{ec2}) gives rise to
finitely many nonrational plane curves, each of which can
have only finitely many integer points.
\eep

The case when $c$ is the successor of a perfect square has
received a lot of attention by people working on Terai's
conjecture (cf.~\cite{dcaod} and the references therein).
Our next result improves on all published results
on this case.

\bet{teunu}
If  in the representation for $c$ derived from Lemma~\ref{lemord}
one has $v=1$, then the system~\emph{(\ref{ec2})} has   no
solutions subject to restrictions from~\emph{(\ref{ec2bis})}.
\eet
\bep
We argue by reduction to absurd. Assume that $c=u^2+1$,
and consequently
$b= \pm  \sum_{j=0}^{(r-1)/2} c_{r,j}\, c^j (-4)^{(r-1)/2-j}$.
Suppose that $(x,y,z)$ is a
solution to the simultaneous equations~(\ref{ec2})
satisfying all the conditions from~(\ref{ec2bis}).
 From Lemma~\ref{lemord} applied for $n=r$ we know that
$a+ib=\eta_1 (u+\eta_2 i)^r$ with $\eta_1$, $\eta_2\in \{\pm1\}$,
thus
\[
a \equiv \pm ru \pmod{u^3},\quad  b \equiv \pm \left(1 - 
\binom{r}{ 2}u^2\right)\pmod{u^4},
\]
and it follows that
\[
c^z=a^2+b^y \equiv r^2 u^2+\left( 1-\frac{1}{2}r(r-1)\, y\, 
u^2\right) \equiv 1+z\, u^2 \pmod{u^4},
\]
that is, $\frac{1}{2}r(r-1)\, y+z\equiv r^2 \pmod{u^2}$.
On noting that the left-hand side of this relation is greater
than the right-hand side (because $y\ge 6$), one obtains the
first inequality from the chain
\beq{ec22}
  u^2+r^2 \le \frac{1}{2}r(r-1)\, y+z < \frac{1}{2}r^2\, y.
\eeq
The second inequality holds  since $z<r\, y/2$.  Indeed, 
$c^{ry/2} =(a^2+b^2)^{y/2}>a^2+b^y=c^z$.
Since in this case $u^2\ge 4\cdot 10^{10}$, Eq.~(\ref{ec22}) 
readily  contradicts the bounds $r< 770$ and $y< 1800$
already obtained.
\eep

We are now in a position to prove Terai's conjecture when
$b$ is a prime power. The proof relies on the observation
that $b$ is of the form $\pm vU_r$, where
\[U_r=U_r(\alpha,\beta)=\frac{\alpha^r-\beta^r}{\alpha-\beta}
\]
is the $r$th Lucas number associated to the pair $(\alpha,
\beta)=(u+v\sqrt{-1},u-v\sqrt{-1})$. In a subsequent proof
we shall use the fact that $a=\pm u\tilde{U}_r$, with
\[\tilde{U}_r=\tilde{U}_r(\tilde{\alpha},\tilde{\beta})=
\frac{\tilde{\alpha}^r-
\tilde{\beta}^r}{\tilde{\alpha}-\tilde{\beta}}
\]
the $r$th Lehmer number associated to the pair $(\tilde{\alpha},
\tilde{\beta})=(u+v\sqrt{-1},-u+v\sqrt{-1})$. Recall that a
prime divisor of $U_r$, respectively $\tilde{U}_r$, is called
primitive if it does not divide
\beq{ec33}
(\alpha-\beta)^2 U_1\cdots U_{r-1}=-4v^2U_1\cdots U_{r-1},
\eeq
respectively
\beq{ec34}
(\tilde{\alpha}^2-\tilde{\beta}^2)^2 \tilde{U}_1\cdots
\tilde{U}_{r-1}=-16u^2v^2\tilde{U}_1\cdots \tilde{U}_{r-1}.
\eeq

Bilu, Hanrot and Voutier~\cite{bhv} showed that for $n>30$,
every $n$th Lucas and Lehmer number has a primitive divisor.
Moreover, they and Abouzaid~\cite{abid} have given the
complete list of $n$ and $(\alpha,\beta)$, respectively
$(\tilde{\alpha},\tilde{\beta})$, for which $U_r(\alpha,\beta)$
or $\tilde{U}_r(\tilde{\alpha},\tilde{\beta})$ does not
have a primitive divisor.

\bet{tecb}
If  $b$  is a prime power then the system~\emph{(\ref{ec2})}
has   no solutions subject to restrictions
from~\emph{(\ref{ec2bis})}.
\eet\bep
Let $p$ be an odd prime and $s$ a positive integer such that
$b=p^s$.  Having in view the result just proved, we conclude
that if the
system~(\ref{ec2}) has a solution satisfying~(\ref{ec2bis}),
then $p$ divides $v$.
Therefore, either $U_r=1$ or its only prime divisor $p$ is not
primitive (see Eq.~(\ref{ec33})). Checking the relevant tables
from~\cite{bhv} and~\cite{abid}, one finds that one necessarily
has $r=3$, $5$, $7$ or $13$. Moreover, when $r=3$, $c$ would
result even, in contradiction to~(\ref{ec2bis}). For $r=5$,
all the candidates for $(\alpha,\beta)$ do not yield an integer
value for $v$, while for $r=7$ or $13$ the resulting value for
$u$ is not integer.
\eep

To the best of our knowledge, the literature contains nothing
of the kind of our next result.

\bet{teca}
If  $a$  is a prime power then the system~\emph{(\ref{ec2})}
has   no solutions subject to restrictions
from~\emph{(\ref{ec2bis})}.
\eet\bep
As explained before, we use the equality $a=\pm u\tilde{U}_r$,
with $u\ge 2$. We proceed as in the previous proof,
reasoning about the Lehmer pair $(\tilde{\alpha}, \tilde{\beta})$
instead of the Lucas pair $(\alpha,\beta)$. Since the
differences are insignificant, the details can be safely left
to the reader.
\eep

\section{Further results}\label{s5}

In subsequent reasonings we shall need to know
that $v_1\ne 1$. This fact follows from the following.

\bel{le55}
With the notation of the previous section we have the two 
following results:
\[
\min\{u_1/v_1,v_1/u_1\}\le 0{.}01 \ \Longrightarrow \ r\le 659 \
\quad   and \quad z\le 845
\]
and
\[
\min\{u_1/v_1,v_1/u_1\}\ge 0{.}001856 .
\]
In particular,
\[
\min\{u_1,v_1\}\ge 372. 
\]
\eel
\bep
With the notation $\eps_1=u_1+iv_1 =|\eps|e^{i\xi_1}$ and  
$\xi'_1=\pi/2-\xi_1$, the corresponding linear form is
\[
\Lambda = z (2i\xi_1) -k (i\pi/2) = z (-2i\xi'_1) -k' (i\pi/2)
\]
and when $\xi_1$ or $\xi'_1$ is small we can get much better 
estimates   in the application of Laurent's lower bound.
Technically: we can take a much larger radius of interpolation 
and we   obtain the above upper bounds for $r$ and~$z$.

The proof of the second result is elementary. We have  
$a=c^{z/2}|\cos(z\xi_1)|=c^{z/2}|\sin(z\xi'_1)|$.
Hence the condition $(z+1)\xi_1<\pi/2$ implies
\[
|\cos(z\xi_1)|\ge \cos(\pi/2-\xi_1)=\sin \xi_1=\frac{v_1}{\sqrt c},
\]
where $0<\xi_1<\tan \xi_1=v_1/u_1$. It follows that
\[
\frac{v_1}{u_1} < \frac{\pi }{2\times 846} =0{.}001856733\ldots\  
\Longrightarrow \ a>c^{z/2-1}\ge c^{r/2}.
\]
Since $a^2+b^2=c^r$, this is a contradiction that proves the 
lower bound   $v_1/u_1\ge 0{.}001856$. A similar reasoning
leads to the inequality  $u_1/v_1\ge 0{.}001856$.

Now, since $u_1^2+v_1^2>4\cdot 10^{10}$, a simple computation 
gives   $\min\{u_1,v_1\}\ge 372$.
\eep

In a similar way we can prove partially analogous results 
concerning the pair $(u,v)$.

\bel{le55b}
The following implication holds
\[
\min\{u/v,v/u\}\le 0{.}01 \ \Longrightarrow \ r\le 553 \quad 
and \quad z\le 705 .
\]

If the Diophantine system~\emph{(\ref{ec2})} has solutions
satisfying~\emph{(\ref{ec2bis})}  with
\[
b\ge c^{(r-1)/2}
\]
(which is true if $v(r+1)<\pi u$) then
\[
y \ge 6 \ \Longrightarrow \ r\le 101\quad  and  \quad z\le 299.
\]
Moreover, again under the hypothesis $b\ge c^{(r-1)/2}$,
\[
y \ge 10 \ \Longrightarrow \ r\le 47\quad and \quad  z\le 227,
\quad
y \ge 14 \ \Longrightarrow \  r\le 31 \quad and \quad  z\le 209,
\]
\[
y \ge 18 \ \Longrightarrow \  r\le 23 \quad and \quad  z\le 197,
\quad
y \ge 22 \ \Longrightarrow \  r\le 19 \quad and \quad  z\le 189,
\]
\[
y \ge 30 \ \Longrightarrow  \ r\le 13  \quad and \quad z\le 185,
\quad
y \ge 50 \ \Longrightarrow  \ r\le \phantom{0}7  \quad and \quad 
z\le 161,
\]
\[
y \ge 70 \ \Longrightarrow  \ r\le \phantom{0}5  \quad and \quad 
z\le 155,
\quad
y \ge 98 \ \Longrightarrow  \ r\le \phantom{0}3  \quad and \quad 
z\le 147,
\]
and there is no solution for $y\ge 142$.
\eel

We add some other estimates related to $b$.

\bel{le55c}
If the Diophantine system~\emph{(\ref{ec2})} has solutions
satisfying~\emph{(\ref{ec2bis})}  then
\[
ry/2=z+2t,\quad  with \quad t\ge 1,
\]
and
\[
b< c^{\frac r 2 -\frac 2 y}.
\]
Moreover, if
\[
b\ge (1+10^{-20}) \, c^{\frac r 2 -\frac 4 y}
\]
then
\[
ry/2=z+2.
\]
When it holds $\,ry/2=z+2\,$  then
\[
y \ge 6\phantom{0} \ \Longrightarrow \ r\le 101,\quad
y \ge 10 \ \Longrightarrow \ r\le 47,\quad
y \ge 14 \ \Longrightarrow \  r\le 29,
\]
\[
y \ge 18 \ \Longrightarrow \  r\le \phantom{0}19, \quad
y \ge 22 \ \Longrightarrow \  r\le 17, \quad
y \ge 26 \ \Longrightarrow \  r\le 13,
\]
\[
y \ge 30 \ \Longrightarrow  \ r\le \phantom{0}11, \quad
y \ge 38 \ \Longrightarrow  \ r\le \phantom{0}9, \quad
y \ge 42 \ \Longrightarrow  \ r\le \phantom{0}7,
\]
\[
y \ge 50 \ \Longrightarrow  \ r\le \phantom{00}5,  \quad
y \ge 66 \ \Longrightarrow  \ r= \phantom{0}3, \quad
\phantom{y \ge 14 \ \Longrightarrow \  r\le 29,}
\]
and there is no solution for $y\ge 102$.
\eel

\bep
We give a proof just for the first two assertions. From the relations
\[
(a^2+b^2)^{y/2}>a^2+b^y=c^z
\]
we deduce $ry/2>z$ and the first assertion follows since $ry/2$ and  
$z$ are both odd.

If $b\ge (1+10^{-20})c^{\frac r 2 -\frac 4 y}$ then, since  
$b^y>(1-10^{-21})c^z$, we see that
$z>ry/2-4$ and the relation $\,z=ry/2-2\,$ follows from the first 
assertion.

The remaining estimates result from computation with the help of 
lower   bounds on linear forms.
\eep

It is very likely that actually there are no solutions
to~(\ref{ec2}) under the conditions stated
in~(\ref{ec2bis}).
This is the case under the  hypothesis of the next result.

\bet{techen}
The system~\emph{(\ref{ec2})} has no solutions $(r,y,z)$
subject to restrictions~\emph{(\ref{ec2bis})} in which $z$
is divisible by $3$ and $y\ne  6$, $10$, $14$, $18$, $30$, $42$, 
$50$, $54$, $62$, $70$, $90$, $98$, $126$,  $150$, $162$, $186$, 
$210$, $250$, $270$, $294$, $310$, $350$, $378$, $434$, $450$, 
$486$, $490$, $558$, $630$.
\eet
\bep
I. Chen~\cite{chen} very recently proved that
for any prime  satisfying the restrictions $7<p <10^7$
and $p\ne 31$ there are no coprime integers $A$, $B$, $C$
satisfying
\[
A^2+B^{2p}=C^{3}.
\]
This confirms Terai's conjecture in case $z$ is
multiple of $3$ and $y$ has  a prime divisor $p>7$, $p\ne 31$.
The only values of the $y$-component in a solution of
Eqs.~(\ref{ec2})--(\ref{ec2bis}) not covered by Chen's result
are listed having in view Proposition~\ref{pr55}.
\eep

\medskip

The following remarks are helpful when trying to further 
reduce the number of candidate pairs $(y,z)$.

\medskip

\noindent
\textbf{Remark 1.}  When $r$ divides $z$, we may remove the
multiples of $3$ from this list because Mignotte and
Peth\H o~\cite{mmp} have proved that if there are points
with  both coordinates greater than 1 on the
curve $X^m-X=Y^n-Y$, then  $\gcd (m,n)=1$.

\medskip

\noindent
\textbf{Remark 2.}
A  deep result of Darmon and M\'erel~\cite{dm}, according to 
which the equation $X^n+Y^n=Z^2$  has no solutions in nonzero
integers when $n\ge 4$,  implies that $\gcd(y,z)\le 3$ always 
holds.

\medskip

Our last result is a bit surprising because it shows that
the hypothesis $b>a$ from the main results of~\cite{ter2},
\cite{dcao} (see $(\alpha)$ and respectively $(\gamma)$ in 
Introduction)  and~\cite{leu} is never fulfilled (the 
reader is warned that  in Le's paper $b$ denotes the unique
even  number among $a$ and $b$).

\bepr{pr55}
If system~\emph{(\ref{ec2})} has  solutions subject to
restrictions  from~\emph{(\ref{ec2bis})} then
\[
a>4.608 \, b , \quad
c>3^{y-10}\quad  and \quad y \le 2z+4.
\]
Moreover
\[
  y \le 2z-4 \quad  for \quad y\ge 34
\]
and
\[
y\ge 602 \ \Longrightarrow \ r\le 149 \quad  and \quad  z \le 319.
\]
In particular, 
\[
 y\le 634.
\]
\eepr
\bep
When $y\le 10$ one has $c>10^y$ because $c>4\cdot 10^{10}$.
For the same reason, $c>2.2^y$ when $y$ is between
14 and 30. It is much harder to obtain similar
inequalities for higher values of $y$. We now prove that
it always  holds $c>2{.}1716^y$.

As seen above, $v_1$ has a prime divisor $p$.
Recall that in Lagrange's formula given in Lemma~\ref{le3}
the coefficients for $n$ odd are
\[
c_{n,j}= \frac{(n-j-1)! \, n}{(n-2j)!\, j!},
\]
where $0\le j \le (n-1)/2$, and the quotient $(n-j-1)! /j!$ is
an integer. It follows that  we have
\[
v_p(c_{n,j}) \ge v_p(n)-v_p \bigl( (n-2j)!\bigr)
> v_p(n) - \frac{n-2j}{p-1} \ge v_p(n) - \frac{n-2j}{2}.
\]
As $p$ divides $v_1$, it  does not divide $c$ and therefore
\begin{align*}
v_p (c_{z,j} c (-4v_1^2)^{(z-1)/ 2-j}) & =  &
v_p (c_{z,j}) +(z-1-2j)v_p(v_1)  \phantom{v_p(z)-(z-1)}\\
& \ge & v_p(z)-\frac{1}{2}(z-1-2j)+(z-1-2j)v_p(v_1) \\
& \ge & v_p(z)+\frac{1}{2}(z-1-2j)v_p(v_1)
> v_p(z) \phantom{v_p(z)-}
\end{align*}
for $0\le j < (z-1)/2$.

Corollary~\ref{cor1} yields
\[
   y\, v_p (b) = 2\bigl( v_p(z)+ v_p(v_1)\bigr).
\]
Having in view the upper bounds for $z$ given in Lemma~\ref{le0tri}, 
we see that for $y\ge 34$ it holds $v_3(z)\le 5$,
and $v_p(z)\le 3$ for $p\ge 5$. Consequently,
for $p=3$ one obtains
\[c>v_1^2 \ge 3^{y-10}\ge  3^{y(1-5/17)}>2{.}1716^y.
\]
For $p\ge 5$ one has $c>p^{y-6}\ge 5^{y-6}>3^{y-10}$, so that
the claim that $c>2{.}1716^y$ is true for any solution of the
system~(\ref{ec2}) satisfying conditions~(\ref{ec2bis}).

Denote provisionally $\mu =b^2/c^r$. Then
$c^z>b^y$  implies that
\[\mu ^{y/2}<c^{z-ry/2}\le c^{-2} <2{.}1716^{-2y}.
\]
Hence,
\[a= b\, \sqrt{\mu ^{-1}-1}>b\, \sqrt{2{.}1716^{4}-1}>4{.}608
\,  b.
\]

Since $a^2<c^{z-2}$ and $c>10^{10}$ we have $b^y>(1-10^{-20})c^z$, 
and   the inequality $c>2{.}1716^y$ implies
\[
b>2{.}171^z.
\]

The inequalities relating $y$ and $z$ are proved in three steps.
First, we show that we always have $y\le 2z+12$. Next, we disprove
the equalities $y=2z+8$ and $y=2z+12$ by combining information
already available with some more computations. Similar arguments
are employed to show that one can not have $y=2z+4$ for $y\ge 34$,
while $y\ne 2z$ follows from the result of Darmon and M\'erel
mentioned in Remark~2. Here are the details.

The upper bound
\[
y < z \left(2+\frac{\log\left(1+(r+1)^2/\pi^2 \right)}{ \log b} 
\right),
\]
combined with the lower bound $b>2{.}171^z$, leads to
\[
y < 2z+\frac{\log\left(1+(r+1)^2/\pi^2 \right)}{ \log 2{.}171}.
\]
The bound $\,y\le 2z+12\,$ is trivially satisfied for $y\le 22$, 
and for $y\ge 22$ we have seen that $r<396$, so that 
\[
y < 2z + \frac{\log \bigl(1+396^2/\pi^2 \bigr)}{\log 2{.}171}
< 2z+12.5,
\]
which implies $\,y\le 2z+12$. To show that the equality in this
relation never holds, one argues similarly to the case  $y=2z+4$ 
for $y\ge  34$ detailed below. Then one repeats the reasoning
to show that $y\ne 2z+8$, so that we always  have
\[
   y \le 2z+4.
\]

Suppose that $y=2z+4$ for some $y\ge34$. Theorem~\ref{techen}
implies $y\ge 38$ and we  verify by
a computation with linear form estimates that
\[
y\ge 38 \ \Longrightarrow \ r\le 239.
\]
 From the relation $b^y<c^z$, we get
\[
y = 2z+4   \ \Longrightarrow \ b < c^{1/2-2/y}.
\]
Besides we know that
\[
b < c^{(r-1)/2} \ \Longrightarrow \  u<(r+1)v/\pi\ \Longrightarrow \   
c <(1+(r+1)^2/\pi^2) v^2.
\]
Put $b=vb'$---then $b'$ is a positive integer. The above facts imply
\[
y = 2z+4 \ \Longrightarrow \   c < \left(\frac{1}{b'}\sqrt{1+
\frac{(r+1)^2}{\pi^2}} \, \right)^{y/2}.
\]
Now we consider $v_1$.
We have $v_1<u_1 \pi /(2(z+1))$, thus 
\[
v_1^2 < \left(1+\frac{\pi^2}{4(z+1)^2} \right)^{-1}c.
\]
Moreover, we can write
\[
v_1 = w_1^{y/2}/w_0, \quad {\rm where}\quad w_0 \mid \gcd(z,v_1^2)
\]
and  $b=w_1 b''$, where $b''$ is a positive integer.

If  $b'=1$, a short computer verification  shows that $w_1<9$ for  
$y\ge 38$. Since $w_1$ is an odd integer greater than $1$, one has
$w_1\in \{3,5,7\}$. But we know that $b$ is not a  power of a prime, 
hence $b'\ge 3$, with $b'\ge 5$ when $w_1=3$. Using now $b'\ge 3$ 
another computer   verification leads to
  $w_1\in \{3,5\}$ and $w_1=3$ for $y\ge 102$.

Now we apply again Laurent's result but with the much better lower  
bound $c>3^{y-10}$
(better for $y\ge 34$ than $c>4\cdot 10^{10}$) and we get for example
\[
y\ge 102 \ \Longrightarrow \ r\le 181 \quad {\rm and}\quad z \le 373,
\]
\[
y\ge 302 \ \Longrightarrow \ r\le 157 \quad {\rm and}\quad z \le 329,
\]
\[
y\ge 602 \ \Longrightarrow \ r\le 149 \quad {\rm and}\quad z \le 319.
\]

Comparing the previous estimates we conclude that we always have 
$y\le 634$.

\eep

\end{document}